\documentclass[fleqn]{mat01}
\usepackage{times,mathtimy,amssymb,boxit}
\usepackage[cmtip,v2]{xy}
\usepackage{xspace}

\begin{document}

\setcounter{page}{139}
\firstpage{139}

\font\xx=msam5 at 10pt
\def\ab{\mbox{\xx{\char'03}}}

\font\sa=tibi at 10.4pt
\def\d{\hbox{d}}

\def\defi{\trivlist\item[\hskip\labelsep{\bf DEFINITION.}]}
\def\remark{\trivlist\item[\hskip\labelsep{\it Remark.}]}
\def\remarks{\trivlist\item[\hskip\labelsep{\it Remarks}]}
\def\noot{\trivlist\item[\hskip\labelsep{{\it Note.}}]}
\def\thoe{\trivlist\item[\hskip\labelsep{{\bf Theorem}}]}

\newtheorem{theo}{Theorem}
\renewcommand\thetheo{\arabic{section}.\arabic{theo}}
\newtheorem{theor}[theo]{\bf Theorem}
\newtheorem{definit}[theo]{\rm DEFINITION}
\newtheorem{lem}[theo]{Lemma}
\newtheorem{propo}[theo]{\rm PROPOSITION}
\newtheorem{rema}[theo]{Remark}
\newtheorem{exam}[theo]{Example}
\newtheorem{coro}[theo]{\rm COROLLARY}
\newcommand{\ch}{\overline{\mathbb H}}
\newcommand{\md}{{\mathrm d}}

\markboth{Pablo Ar\'es-Gastesi and Indranil Biswas}{Jacobian of a Klein surface}

\title{The Jacobian of a nonorientable Klein surface}

\author{PABLO AR\'ES-GASTESI and INDRANIL BISWAS}

\address{School of Mathematics, Tata Institute of
Fundamental Research, Mumbai 400 005, India\\
\noindent E-mail: pablo@math.tifr.res.in; indranil@math.tifr.res.in}


\volume{113}

\mon{May}

\parts{2}

\Date{MS received 10 July 2002; revised 17 January 2003}

\begin{abstract}
Using divisors, an analog of the Jacobian for a compact connected
nonorientable Klein surface $Y$ is constructed. The Jacobian is
identified with the dual of the space of all harmonic real one-forms on
$Y$ quotiented by the torsion-free part of the first integral homology
of $Y$. Denote by $X$ the double cover of $Y$ given by orientation. The
Jacobian of $Y$ is identified with the space of all degree zero
holomorphic line bundles $L$ over $X$ with the property that $L$ is
isomorphic to $\sigma^*\overline{L}$, where $\sigma$ is the involution
of $X$.
\end{abstract}


\maketitle

\section{Introduction}

Let $Y$ be a compact connected nonorientable Riemann surface,
that is, each transition function is either holomorphic or
anti-holomorphic. We consider surfaces without boundary.
Let $X$ denote the double cover of $Y$ given
by the local orientations. So $X$ is a compact connected
Riemann surface.

In \S\ref{sect.-1}, we define a morphism from $Y$ to $\ch$, the closure
of the upper half-plane in the Riemann sphere $\widehat{\mathbb C}$.
Let $\text{Div}_0(Y)$ denote the group defined by all formal finite sums of
the form $\sum n_iy_i$, where $n_i\in {\mathbb Z}$ with $\sum n_i=0$ and
$y_i\in Y$. We call such a divisor $D$ to be principal if there is a
morphism (see \S\ref{sect.-1} for the definition of morphism) $u$ from $Y$ to $\ch$ with the property that
\begin{equation*}
D\,=\, u^{-1}(0) - u^{-1}(\infty)\,.
\end{equation*}
Let $J_0(Y)$ denote the quotient of $\text{Div}_0(Y)$ by its subgroup
consisting of all principal divisors. This $J_0(Y)$ is the analog of the
Jacobian for a nonorientable Riemann surface.

Harmonic one-forms are defined on $Y$. Let ${H}^1_{\mathbb
R}(Y)$ denote the space of all harmonic real one-forms on $Y$. The
torsion-free part of ${H}_1(Y, {\mathbb Z})$ is a subgroup of
${\mathcal H}^1_{\mathbb R}(Y)^*$. The quotient is identified with
$J_0(Y)$. This is proved by showing that ${\mathcal H}^1_{\mathbb R}(Y)$
is identified with the space of all holomorphic one-forms $\omega$ on
$X$ satisfying the identity $\overline{\omega}\,=\, \sigma^*\omega$,
where $\sigma$ is the nontrivial automorphism of the double cover $X$ of
$Y$ (Theorem~\ref{thm:harmonic}).

For a holomorphic line bundle $L$ over $X$, the pullback $\sigma^*
\overline{L}$ is again a holomorphic line bundle over $X$. We show that
$J_0(Y)$ is identified with the group of all holomorphic line bundles
$L$ over $X$ for which the holomorphic line bundle $\sigma^*
\overline{L}$ is isomorphic to $L$ (Theorem~\ref{thm.-4.1}).

A compact Riemann surface is a smooth projective curve over $\mathbb C$.
Conversly, every smooth projective curve over $\mathbb C$ corresponds to
a compact Riemann surface. If we take a smooth projective curve
$X_{\mathbb R}$ defined over $\mathbb R$, then using the inclusion of
$\mathbb R$ in $\mathbb C$ we get a smooth projective curve $X_{\mathbb
C}$ over $\mathbb C$. Now, since the involution of $\mathbb C$ defined
by conjugation fixes $\mathbb R$, the complex curve $X_{\mathbb C}$ is
equipped with an anti-holomorphic involution that reverses the
orientation. Conversely, every complex projective curve equipped with an
anti-holomorphic involution is actually defined over $\mathbb R$. If
the involution does not have any fixed points, that is, the curve does
not have any real points, then it is called an imaginary curve.

Therefore, a nonorientable Riemann surface $Y$ (without boundary) corresponds to an
imaginary algebraic curve defined over $\mathbb R$. The Jacobian of the
complexification $Y{\mathbb C}$ is also the complexification of a
variety defined over $\mathbb R$. The Jacobian $J_0(Y)$ coincides with
this variety defined over $\mathbb R$.

\section{Divisors on a nonorientable surface}\label{sect.-1}
\setcounter{theo}{0}
Let $Y$ be a compact connected nonorientable surface. In other words,
$Y$ is a compact connected nonorientable smooth manifold of dimension
two, and $Y$ has a covering by smooth coordinate charts such that each
transition function is either holomorphic or anti-holomorphic. Any
coordinate chart in the maximal atlas satisfying the above condition on
transition functions will be called \textit{compatible}. Such a
nonorientable surface is called a \textit{Klein surface}.

\begin{definit}\label{defn:divisor}$\left.\right.$\vspace{.5pc}

\noindent {\rm A \textit{divisor} $D$ on $Y$ is a formal sum of type
\begin{equation*}
D \,=\, \sum_{y\in Y} n_yy\, ,
\end{equation*}
where $n_y\in \mathbb Z$ and $n_y=0$ except for a finitely
many points of $Y$.}
\end{definit}

\begin{definit}$\left.\right.$\vspace{.5pc}

\noindent{\rm The \textit{degree} of a divisor $D = \sum_{y\in Y} n_y y$
is defined to be the integer ${\rm deg}(D) \,:=\,
\sum_{y\in Y} n_y$.}
\end{definit}

We will denote by $\text{Div}(Y)$ the set of all divisors on $Y$. Let
$\text{Div}_d(Y)\, \subset\, \text{Div}(Y)$ be the divisors of degree
$d$.

Let $\pi:X\to Y$ be a double cover of $Y$ given by local orientations on
$Y$. So for a contractible open subset $U\, \subset\, Y$, the inverse
image $\pi^{-1}(U)$ is two copies of $U$ with the two possible
orientations on $U$ (see \cite{ares:klein} for more details on Klein
surfaces and their double covers).

Therefore, $X$ is a Riemann surface, and the change of orientation
defines an anti-holomorphic involution $\sigma:X\to X$ that commutes
with $\pi$.

The involution $\sigma$ induces in a natural way a mapping on
the set of divisors on the Riemann surface $X$ as follows
\begin{align*}
\sigma^*\,:\,\text{Div}(X) & \longrightarrow \text{Div}(X) \\
\sum m_j\, x_j & \longmapsto \sum m_j\, \sigma(x_j).
\end{align*}
Observe that $\sigma^*$ preserves the degree.

Similarly, the quotient map $\pi:X\to Y$ induces mappings between the
divisors on $X$ and $Y$. To define those mappings we first set up some
notation. For any
point $y \in Y$ we will denote by $\pi^{-1}(y)$ the divisor
given by the inverse image of $y$.
In other words, $\pi^{-1}(y) = x + \sigma(x)$, where
$x\in X$ is a point satisfying $\pi(x) = y$. Then we can define two
mappings as follows:
\begin{align*}
\pi^*: \text{Div}(Y) & \to \text{Div}(X), \qquad \pi_*: \text{Div}(X)  \to
\text{Div}(Y),\\
\sum_{j=1}^s n_j\,y_j & \mapsto \sum_{j=1}^s n_j\,\pi^{-1}(y_j), \qquad
\sum_{j=1}^s m_j\,x_j  \mapsto \sum_{j=1}^s m_j\,\pi(x_j).
\end{align*}
Observe that $(\pi_*\circ\pi^*)(D) = 2D$ and $(\pi^*\circ\pi_*)(E) =
E+\sigma^*(E)$ for $D\in \text{Div}(Y)$ and $E \in \text{Div}(X)$.

Let $\text{Div}(X)^{\sigma^*}$ denote the set of fixed points of
$\sigma^*$ on $\text{Div}(X)$.

The following lemma follows immediately from the above definitions.

\begin{lem}\label{lemma-2.1}
The group ${\rm Div}(Y)$ is identified with
${\rm Div}(X)^{\sigma^*}$. The isomorphism takes the subgroup
${\rm Div}_0(Y)$ to ${\rm Div}(X)_0^{\sigma^*}
= {\rm Div}_0(X)\cap {\rm Div}(X)^{\sigma^*}$.
\end{lem}

Let $j:\widehat{\mathbb C} \to \widehat{\mathbb C}$ denote the mapping
induced by conjugation on the Riemann sphere $\widehat{\mathbb C}$, so
that $j(z) = \bar{z}$ and $j(\infty) = \infty$. The quotient space is a
surface with boundary, $\overline{\mathbb H} = \widehat{\mathbb C}/
\langle j\rangle $. We can also identify $\overline{\mathbb H}$ with
the closure of $\mathbb H$ (the upper half-plane)
in the Riemann sphere. Let
\begin{equation*}
p\,:\,\widehat{\mathbb C} \,\longrightarrow\,
\overline{\mathbb H}
\end{equation*}
denote the quotient map. After identifying $\overline{\mathbb H}$ with
the closure of $\mathbb H$ the map $p$ coincides with the
one defined by $p(x+\sqrt{-1}y) = x +\sqrt{-1}|y|$ and
$p(\infty) =\infty$.

A \textit{morphism} from $Y$ to $\ch$ is a continuous mapping
\begin{equation*}
u\,:\,Y\, \longrightarrow\, \ch
\end{equation*}
such that if $(U,w)$ is a local coordinate function
defined on $Y$, compatible with the Riemann surface
structure, with $w(U) \subset \mathbb H$,
then there exists a holomorphic function $F:w(U) \to
\mathbb C$ that makes the following diagram commutative:
\begin{equation*}
\begin{split}
\xymatrix{
U \ar[r]^u \ar[d]_w & \ch \\
{\mathbb H} \ar[r]_F & {\mathbb C} \ar[u]_p
}
\end{split},
\end{equation*}
where $p$ is defined above.

Let $u$ be a morphism, as above, from $Y$ to $\ch$ which is not
identically equal to $0$ or $\infty$. If $z_0$ is a point of $\ch$, then
by $u^{-1}(z_0)$ we understand the divisor given by the inverse image of
$z_0$ under $u$ (so the integers $n_j$ in Definition~\ref{defn:divisor}
are given by the multiplicities of $u$ at the corresponding points).
Since $0$ and $\infty$ in $\widehat{\mathbb C}$ project to two different
points on $\ch$,
\begin{equation*}
\text{div}(u)\, :=\, u^{-1}(0) - u^{-1}(\infty)\,\in\,
\text{Div}(Y)
\end{equation*}
is a divisor on $Y$.

\begin{definit}$\left.\right.$\vspace{.5pc}

\noindent {\rm A divisor $D \in {\rm Div}(Y)$ is called {\it principal} if
$D = \text{div}(u)$ for some morphism
$u:Y \to \ch$ of the above type. The set of principal divisors of
$Y$ will be denoted by $\text{Div}_P(Y)$.}
\end{definit}

\begin{propo}\label{prop.-2.1}$\left.\right.$\vspace{.5pc}

\noindent A divisor $D$ on $Y$ is principal if and only if there
exists a divisor
$E \in {\rm Div}_P(X) \cap {\rm Div}(X)^{\sigma^*}$ with
$\pi^*D = E$.
\end{propo}

\begin{proof}
Let $E\, =\, \text{div}(f)$ be a principal divisor in
$\text{Div}(X)^{\sigma^*}$, where $f$ is a non-constant
meromorphic function on $X$. Consider the function
$\psi$ on $X$ defined by
$\psi(x)\, =\, \overline{f(\sigma(x))}$
on $X$. This function $\psi$ is clearly meromorphic.

Since $E\in \text{Div}(X)^{\sigma^*}$, we have
$\text{div}(\psi)\, =\,\text{div}(f)$. Consequently,
there exists a constant $c\in {\mathbb C}\setminus \{0\}$
such that $\psi \, =\, cf$.

Therefore, we have $f(x) = \psi(x)/c = \overline{f(\sigma(x))}/c
= \overline{\psi(\sigma(x))}/|c|^2 = f(x)/|c|^2 $. Take
$c_0\in {\mathbb C}$ with $c^2_0 = c$. Set
$f_0 \, =\, c_0f$.

The divisor for the meromorphic function $f_0$ coincides
with $E$. Furthermore, $f_0$ satisfies the condition
\begin{equation*}
f_0\circ\sigma \, =\, \overline{f_0}\, .
\end{equation*}
Therefore, it induces a map
\begin{equation*}
\hat{f}\, :\, Y\, :=\, X/\sigma \, \longrightarrow\,
\overline{\mathbb H}\, :=\, \widehat{\mathbb C}/\langle j\rangle
\end{equation*}
with $\text{div}(\hat{f}) = D$.

Conversely, let $D\, =\, \text{div}(u)$ be a principal divisor on $Y$.
Consider the composition $u\circ\pi\, :\, X\, \longrightarrow\,
\overline{\mathbb H}$. It is straight-forward to see that the function
$u\circ p$ lifts to a smooth function
\begin{equation*}
f\, :\, X\, \longrightarrow\, \widehat{\mathbb C}
\end{equation*}
such that $p\circ f\, =\, u\circ\pi$. There are two such smooth lifts;
one is holomorphic and the other is anti-holomorphic ($u\circ p$ also
has a continuous lift, defined by the inclusion of $\overline{\mathbb
H}$ in $\widehat{\mathbb C}$ which is not smooth). Let $f$ denote the
holomorphic one. Since $\text{div}(f)\, =\, \pi^*(D)\,\in\,
\text{Div}(X)^{\sigma^*}$, the proof of the proposition is
complete.\hfill \ab
\end{proof}

\begin{definit}\label{def.-Jac.-Y}$\left.\right.$\vspace{.5pc}

\noindent {\rm The quotient of $\text{Div}_0(Y)$, the group of
all degree zero divisors on $Y$, by the subgroup of all
principal divisors on $Y$ is called the \textit{Jacobian}
of $Y$. The Jacobian of $Y$ will be denoted by $J_0(Y)$.}
\end{definit}

From Proposition~\ref{prop.-2.1}, it follows immediately that by sending
any divisor $D$ on $Y$ to the divisor $\pi^*D$ on $X$ we obtain an
injective homomorphism from $J_0(Y)$ to the Jacobian $J_0(X)$ of $X$.
From Lemma~\ref{lemma-2.1}, it follows that $J_0(Y)$ coincides with the
fixed point set of the involution of $J_0(X)$ defined by $\sigma$.

A function $f:W\to \mathbb R$, defined on an open subset of $Y$
is called \textit{harmonic} if for every point $y \in W$, there
exists a compatible coordinate chart $(U,w)$, with
\begin{equation*}
y\,\in\, U\,\subseteq\, W\,,
\end{equation*}
such that the function $f\circ w^{-1}$ is harmonic. Since precomposition
with holomorphic and anti-holomorphic functions preserve harmonicity,
we conclude that harmonic functions are well-defined on $Y$.

We say that a real one-form $\eta$ on $Y$ is \textit{harmonic} if it is
locally given by $\md f$, where $f$ is a harmonic function.

Let $\Omega$ denote the holomorphic cotangent bundle of the Riemann
surface $X$. If $\omega \in {H}^0(X,\Omega)$ is given locally by
$\omega = f \md z$, where $f$ is a holomorphic function, then define
\begin{equation*}
\overline{\sigma^*\omega}
\,:=\, (\overline{f\circ\sigma})\:\md (\overline{z}\circ\sigma)\, .
\end{equation*}
So if $\omega$ is defined over $U$, then $\overline{\sigma^*\omega}$ is
a holomorphic one-form defined over $\sigma (U)$. More generally, for a
one-form $\alpha = u\,\md z+ v\,\md \bar{z}$, set
\begin{equation*}
\sigma^*\alpha = (u\circ\sigma)\:\md (z\circ\sigma) +
(v\circ\sigma)\:\md (\overline{z}\circ\sigma)\, .
\end{equation*}

Let ${\mathcal H}^1_{\mathbb R}(Y)$ and ${\mathcal H}^1_{\mathbb R}(X)$
denote the space of all real harmonic one-forms on $Y$ and $X$
respectively. Using the map $\pi:X\to Y$, we can lift harmonic forms on
$Y$ to smooth forms on $X$. It is easy to see that the pullback of a
harmonic form on $Y$ is a harmonic form on $X$. Therefore, there is a
well-defined injective homomorphism $\pi^*\,:\,{\mathcal H}^1_{\mathbb
R}(Y) \, \longrightarrow\, {\mathcal H}^1_{\mathbb R}(X)$.

The complex structure on $X$ defines a Hodge-$*$ operator on
one-forms on $X$. In local holomorphic coordinates the Hodge-$*$ operator is
\begin{equation*}
*(u\,\md z+
v\,\md\bar{z}) = -\sqrt{-1} u\,\md z + \sqrt{-1} v\md\bar{z}
\end{equation*}
or $*(a\,\md x + b\,\md y) = -b\,\md x + a\,\md y$.

A holomorphic one-form $\omega$ on $X$ will be called
$\sigma$-\textit{invariant} if $\sigma^*\omega \,=\,
\overline{\omega}$. The space of all $\sigma$-invariant forms on $X$
will be denoted by ${H}^0(X,\Omega)^{\overline{\sigma^*}}$.

\begin{theor}[\!]\label{thm:harmonic}
A holomorphic form $\omega \in {H}^0(X,\Omega)$ is
$\sigma$-invariant if and only if there exists a form $\eta \in
{\mathcal H}^1_{\mathbb R}(Y)$ such that $\omega = \beta +
\sqrt{-1}(*\beta)${\rm ,} where $\beta = \pi^*\eta$.

The homomorphism ${H}^1_{\mathbb R}(Y)\, \longrightarrow\,
{H}^0(X,\Omega)^{\overline{\sigma^*}}$ defined by
\begin{equation*}
\eta\, \longmapsto\, \pi^*\eta+ \sqrt{-1}(*\pi^*\eta)
\end{equation*}
is an isomorphism of real vector spaces.
\end{theor}

\begin{proof}
Take any $\omega \in {H}^0(X,\Omega)$. Let $\omega = \beta +
\sqrt{-1}(*\beta)$, where $\beta$ is a real one-form. Now the condition
$\sigma^*\omega \,=\, \overline{\omega}$ immediately implies that
$\sigma^*\beta \,=\, \beta$. Therefore, $\beta$ is the pullback of a
form on $Y$. For any $\eta\in {\mathcal H}^1_{\mathbb R}(Y)$, the form
$\pi^*\eta+ \sqrt{-1}(*\pi^*\eta)$ is a $\sigma$-invariant holomorphic
one-form.

Let
\begin{equation*}
\varphi\, :\, {\mathcal H}^1_{\mathbb R}(Y)\, \longrightarrow\,
{H}^0(X,\Omega)^{\overline{\sigma^*}}
\end{equation*}
be the homomorphism that sends any harmonic form $\eta\in {\mathcal
H}^1_{\mathbb R}(Y)$ to the holomorphic form $\pi^*\eta+
\sqrt{-1}(*\pi^*\eta)$. This homomorphism is injective since a
holomorphic one-form with vanishing real part must be identically zero.

The inverse homomorphism
\begin{equation*}
{H}^0(X,\Omega)^{\overline{\sigma^*}}
\, \longrightarrow\, {\mathcal H}^1_{\mathbb R}(Y)
\end{equation*}
sends a $\sigma$-invariant form $\omega$ on $Y$ to $\eta$ with the
property
\begin{equation*}
\pi^*\eta \, =\, \frac{\omega+\overline{\omega}}{2}\, .
\end{equation*}
This completes the proof of the theorem.\hfill \ab
\end{proof}

\section{The Jacobian}\label{sect.-2}
\setcounter{theo}{0}
A closed oriented smooth path $\gamma$ on $X$ gives an element
$L_\gamma\,\in\, {H}^0(X,\Omega)^*$ defined by
\begin{equation*}
L_\gamma (\omega) \,=\, \int_\gamma \omega,
\end{equation*}
where $\omega\in {H}^0(X,\Omega)$. Using Stokes' theorem we get a
mapping from ${H}_1(X,{\mathbb Z})$ to ${ H}^0(X,\Omega)^*$. The
quotient space ${H}^0(X,\Omega)^* / {H}_1(X,{\mathbb Z})$ will be
denoted by $J_1(X)$.

As we saw in the previous section, for a holomorphic one-form $\omega$
on $X$, the form $\overline{\sigma^*\omega}$ is again a holomorphic
one-form. This involution of ${H}^0(X,\Omega)$ induces an
involution
\begin{equation*}
\sigma_1\, :\, {H}^0(X,\Omega)^*\, \longrightarrow\,
{H}^0(X,\Omega)^*\, .
\end{equation*}
In other words, $(\sigma_1(L))(\omega) =
\overline{L(\overline{\sigma^*(\omega)})}$. It is easy to check that for
any closed smooth-oriented path $\gamma$ on $X$, the identity
\begin{equation*}
\sigma_1(L_\gamma) = L_{\sigma(\gamma)}
\end{equation*}
is valid. So, the involution $\sigma_1$ preserves the subgroup
${H}_1(X,{\mathbb Z})\subset {H}^0(X,\Omega)^*$.

Consequently, the involution $\sigma_1$ of ${H}^0(X,\Omega)^*$
induces an involution on the quotient space $J_1(X)$. The involution of
$J_1(X)$ obtained this way will also be denoted by $\sigma_1$.

Let $g$ be the genus of the compact connected Riemann surface $X$.
Suppose we have a canonical basis of ${H}_1(X,{\mathbb Z})$, say
$\{\alpha_1,\ldots,\alpha_g, \beta_1,\ldots,\beta_g\}$. This means that
the corresponding intersection matrix is
\begin{equation*}
J = \left(\begin{array}{@{}cc@{}}
0 & -{I} \\
{I} & 0
\end{array}\right),
\end{equation*}
where ${I}$ is the identity matrix of rank $g$. Then there exists a
unique basis of ${H}^0(X,\Omega)$, say
$\{\omega_1,\ldots ,\omega_g\}$, such that $\int_{\alpha_k} \omega_j =
\delta_{jk}$ (\cite{fk:book}, Proposition III.2.8). We say that this
basis is \textit{adapted} to the given basis of homology.

Using this adapted basis we
can identify ${H}^0(X,\Omega)^*$ with ${\mathbb C}^g$ by sending
the element $L$ of ${H}^0(X,\Omega)^*$ to the vector
$(L(\omega_1),\ldots ,L(\omega_g))$.

Therefore, for any $\gamma\,\in\, {H}_1(X,
{\mathbb Z})$, we may identify the
element $L_\gamma\in {H}^0(X,\Omega)^*$ with
\begin{equation*}
(L_\gamma(\omega_1),\ldots ,L_\gamma(\omega_g))
\,\in\, {\mathbb C}^g\, .
\end{equation*}
Denote by $\mathcal L$ the lattice in ${\mathbb C}^g$ defined
by ${H}_1(X, {\mathbb Z})$ using this identification.
The quotient space $J_1(X)$ defined earlier is clearly
identified with the quotient ${\mathbb C}^g/{\mathcal L}$.

Assume that the basis $\{\omega_j\}$ is $\sigma$-invariant, that is,
$\overline{\sigma^*(\omega_j)} \,=\, \omega_j$ for each $j\in [1\, ,g]$.
It is easy to check that by the above isomorphism of ${H}^0(X,\Omega)^*$
with ${\mathbb C}^g$ the involution $\sigma_1$ of ${H}^0(X,\Omega)^*$
(defined earlier) coincides with the conjugation defined as $(z_1,
\ldots ,z_g) \, \longmapsto\, (\overline{z_1},\ldots ,\overline{z_g})$.

We will denote by $\sigma_\#$ the involution of ${H}_1(X,
{\mathbb Z})$ induced by the involution $\sigma$ of $X$. Let
\begin{equation*}
\{\gamma_1,\ldots,\gamma_g, \delta_1,\ldots,\delta_g\}
\end{equation*}
be a canonical basis of ${H}_1(X, {\mathbb Z})$
satisfying the condition
$\sigma_\#(\gamma_j) = \gamma_j$ for all $j\,\in\, [1\, ,g]$.
Let $\{\omega_1,\ldots, \omega_g\}$ denote the
corresponding adapted basis.

\setcounter{theo}{0}
\begin{propo}\label{prop.-3.1}$\left.\right.$\vspace{.5pc}

\noindent The above adapted basis $\{\omega_1,\ldots,\omega_g\}$ is
$\sigma$-invariant.
\end{propo}

\begin{proof}
Since
\begin{equation*}
\int_{\sigma_\#\gamma} \omega\, =\, \int_\gamma \sigma^*\omega
\,=\, \overline{\int_\gamma \sigma^*\overline{\omega}}
\end{equation*}
(as $\sigma$ is an involution), the proposition follows
immediately.\hfill \ab
\end{proof}

As in \S\ref{sect.-1}, let $J_0(X)$ denote the quotient
$\text{Div}_0(X)/\text{Div}_P(X)$. For a meromorphic
function $f$ we have
$\sigma^*(\text{div}(f)) = \text{div}(\overline{f\circ\sigma})$.
So $\sigma^*$ induces an involution on $J_0(X)$. This
involution of $J_0(X)$ will be denoted by $\sigma_0$.

Let $\{\omega_1\ldots,\omega_g\}$ be the basis in
Proposition~\ref{prop.-3.1}. Recall the quotient $J_1(X)$ of ${
H}^0(X,\Omega)^*$ defined earlier. The Abel--Jacobi map $A:X \to J_1(X)$
is defined as follows: choose a point $x_0$ of $X$ and set $A(x) =
\big[\int_{x_0}^x\omega_1, \ldots, \int_{x_0}^x\omega_g\big]$, where
the brackets denote the equivalence class in $J_1(X)$. We have
\begin{align*}
A(\sigma(x)) &= \left[ \int_{x_0}^{\sigma(x)} \omega_1, \ldots,
\int_{x_0}^{\sigma(x)} \omega_g \right]
 = \left[ \int_{x_0}^{\sigma(x_0)} \omega_1, \ldots,
\int_{x_0}^{\sigma(x_0)} \omega_g \right]\\
&\quad \ +
\left[ \int_{\sigma(x_0)}^{\sigma(x)} \omega_1, \ldots,
\int_{\sigma(x_0)}^{\sigma(x)} \omega_g \right]\\
& =
c_0 + \left[ \int_{\sigma(x_0)}^{\sigma(x)} \sigma^*(\omega_1), \ldots,
\int_{\sigma(x_0)}^{\sigma(x)} \sigma^*(\omega_1) \right]\\
& =
c_0 + \left[ \int_{x_0}^x \overline{\omega_1}, \ldots, \int_{x_0}^x
\overline{\omega_g} \right] = c_0 + \overline{A(x)},
\end{align*}
where $c_0 = A(\sigma(x_0))$. For a divisor $D = \sum_{j=1}^r n_j\,x_j$,
we define
\begin{equation*}
A(D) \,=\, \sum_{j=1}^r n_j\,A(x_j)\, .
\end{equation*}
If $D$ has degree equal to $0$ then we can write it as $D = \sum_{j=1}^s
x_j - \sum_{j=1}^s y_j$, where $x_j \neq y_k$ (though we can have
repetitions among the $x_j$s or the $y_k$s). Then it is easy to check
that
\begin{equation}\label{eq:abel}
A(\sigma_0 (D)) = \overline{A(D)} = \sigma_1(A(D))\,,
\end{equation}
where $\sigma_1$ and $\sigma_0$ are the earlier defined involutions of
${H}^0(X,\Omega)^*$ and $J_0(X)$ respectively.

By Abel's theorem, the map $A$ can be extended to a map from $J_0(X)$ to
$J_1(X)$. By the Abel--Jacobi inversion problem, the map $A:J_0(X) \to
J_1(X)$ is surjective. Thus \eqref{eq:abel} says that $\sigma_0$ and
$\sigma_1$ are equivalent under $A$, that is, the following diagram
commutes:
\begin{equation}\label{diag.-comm.}
\begin{split}
\xymatrix{
J_0(X) \ar[d]_{\sigma_0} \ar[r]^A & J_1(X) \ar[d]^{\sigma_1} \\
J_0(X) \ar[r]_A & J_1(X).
}
\end{split}
\end{equation}

In the paragraph following Definition~\ref{def.-Jac.-Y} we noted that
the Jacobian $J_0(Y)$ coincides with the fixed point set of $J_0(X)$
for the action of the involution $\sigma_0$. Let $J_1(X)^{\sigma_1}\,
\subset\, J_1(X)$ be the fixed point set for the action of the
involution $\sigma_1$ on $J_1(X)$. From the commutativity of the diagram
in \eqref{diag.-comm.} it follows immediately that $J_0(Y)$ is
identified with $J_1(X)^{\sigma_1}$. Finally using
Theorem~\ref{thm:harmonic}, the Jacobian $J_0(Y)$ is identified with the
quotient of ${\mathcal H}^1_{\mathbb R}(Y)$ by the torsion-free part of
${H}_1(Y, {\mathbb Z})$.

\section{Line bundles on a Klein surface}
\setcounter{theo}{0}
Let $L$ be a holomorphic line bundle over a Riemann surface $X$. By
$\overline{L}$ we will mean the $C^\infty$ complex line bundle over $X$
whose transition functions are the conjugations of the transition
functions for $L$. To explain this, let $U_i$, $i\in I$, be an open
covering of $X$ and assume that over each $U_i$ we are given a
holomorphic trivialization of $L$. So for any ordered pair $i,j\in I$,
we have the corresponding transition\break function
\begin{equation*}
f_{i,j}\, :\, U_i\cap U_j \, \longrightarrow\, {\mathbb C}^*
\end{equation*}
which is holomorphic. The $C^\infty$ complex line bundle $\overline{L}$
has $C^{\infty}$ trivializations over each $U_i$, $i\in I$, and for any
ordered pair $i,j\in I$ the corresponding transition function is
$\overline{f_{i,j}}$. It is easy to see that the collection
$\{\overline{f_{i,j}}\}_{i,j\in I}$ satisfy the cocycle condition to
define a $C^\infty$ complex line bundle.

The line bundle $\overline{L}$ can also be defined without using local
trivializations. A $C^\infty$ complex line bundle is a $C^\infty$ real
vector bundle of rank two together with a smoothly varying complex
structure on the fibers (which are real vector spaces of dimension two).
The underlying real vector bundle of rank two for $\overline{L}$
coincides with the one for $L$. For any $x\in X$, if $J_x$ is the
complex structure on the fiber $L_x$, then the complex structure of the
fiber $\overline{L}_x$ is $-J_x$.

As in \S\ref{sect.-1}, let $Y$ be a nonorientable Klein surface and $X$
its double cover, which is a connected Riemann surface of genus $g$.

Let $L$ be a holomorphic line bundle over $X$. The complex line bundle
$\sigma^*\overline{L}$ has a natural holomorphic structure, where
$\sigma$, as before, is the involution of $X$. To construct the
holomorphic structure on $\sigma^*\overline{L}$, observe that if $f$ is
a holomorphic function on an open subset $U$ of $X$, then
$\overline{f\circ \sigma}$ is a holomorphic function of $\sigma (U)$. We
can choose the above open subsets $U_i$ (sets over which $L$ is
trivialized) in such a way that $\sigma (U_i) \,=\, U_i$. Now, since
each $\overline{f_{i,j}\circ \sigma}$ is a holomorphic function on
$U_i\cap U_j$, the complex line bundle $\sigma^*\overline{L}$ gets
equipped with a holomorphic structure.

\begin{propo}\label{prop.-4.1}$\left.\right.$\vspace{.5pc}

\noindent Let $D$ be a divisor on $X$ of degree $d$ and $L$ the
corresponding holomorphic line bundle ${\mathcal O}_X(D)$
over $X$ of degree $d$. Then the holomorphic line bundle
$\sigma^*\overline{L}$ corresponds to the divisor
$\sigma(D)${\rm ,} that is{\rm ,} $\sigma^*\overline{L}\, \cong\,
{\mathcal O}_X(\sigma(D))$.
\end{propo}

\begin{proof}

Since $L\, \cong\,{\mathcal O}_X(D)$, we have a meromorphic section $s$
of $L$ with the positive part of $D$ as the zeros of $s$ (of order given
by multiplicity) and the negative part of $D$ as the poles of $s$ (of
order given by multiplicity). Since $L$ and $\overline{L}$ are
identified as real rank two vector bundles, the pullback $\sigma^*s$
defines a smooth section of $\sigma^*\overline{L}$ over the complement
(in $X$) of the support of $D$.

It is straight-forward to check that the section $\sigma^*s$ of
$\sigma^*\overline{L}$ is meromorphic. The divisor defined by the
meromorphic section $\sigma^*s$ clearly coincides with $\sigma (D)$.
Consequently, $\sigma^*\overline{L}$ is holomorphically isomorphic to
the line bundle over $X$ defined by the divisor $\sigma (D)$. This
completes the proof of the proposition.\hfill \ab
\end{proof}

Recall the quotient space $J_0(X)\, :=\, \text{Div}_0(X)/
\text{Div}_P(X)$ considered in \S\ref{sect.-1}. The Jacobian
$J_0(X)$ is identified with the space of all isomorphism classes of
degree zero holomorphic line bundles over $X$. The isomorphism sends any
divisor $D$ to the line bundle ${\mathcal O}_X(D)$. As in
\S\ref{sect.-2}, let $\sigma_0$ denote the involution of $J_0(X)$
defined by $\sigma$. From Proposition~\ref{prop.-4.1}, it follows
immediately that the above identification of $J_0(X)$ with degree zero
line bundles takes the involution $\sigma_0$ to the involution defined
by $L\, \longmapsto\, \sigma^*\overline{L}$ on the space of all
isomorphism classes of degree zero line bundles.

Let $D$ be a divisor of degree zero on the nonorientable Klein surface
$Y$. From Proposition~\ref{prop.-2.1}, it follows immediately that $D$ is
principal if and only if $\pi^*D$ is principal. Therefore, we have an
injective homomorphism
\begin{equation}\label{hom.-jac.}
\rho \, :\, \frac{\text{Div}_0(Y)}{\text{Div}_P(Y)}\,
\longrightarrow\, \frac{\text{Div}_0(X)}{\text{Div}_P(X)}
\, =\, J_0(X)
\end{equation}
defined by $D\, \longmapsto\, \pi^*D$, where $\text{Div}_P(Y)$ denotes
the group of principal divisors on $Y$ (as before, $\text{Div}_0$
denotes degree zero divisors).

\begin{theor}[\!]\label{thm.-4.1}
The image of the homomorphism $\rho$ in \eqref{hom.-jac.}
coincides with the subgroup of $J_0(X)$
defined by all holomorphic line
bundle $L$ with $\sigma^*\overline{L}$
holomorphically isomorphic to~$L$.
\end{theor}

\begin{proof}
Let $D$ be a divisor on $Y$ of degree zero. The divisor $\pi^*D$ on $X$
is left invariant by the action of the involution $\sigma$. From the
above remark that the involution $\sigma_0$ is taken into the involution
defined by $L\, \longmapsto\, \sigma^*\overline{L}$, it follows
immediately that the holomorphic line bundle $L\, =\, {\mathcal
O}_X(\pi^*D)$ over $X$ corresponding to the divisor $\pi^*D$ satisfies
the condition $L\, \cong\,\sigma^*\overline{L}$.

For the converse direction, take a holomorphic line bundle $L$ over $X$
which has the property that $\sigma^*\overline{L}$ is isomorphic to $L$.
Let $s$ be a nonzero meromorphic section of $L$. If the divisor
$\text{div}(s)$ is left invariant by the involution $\sigma$, then $L$
is in the image of $\rho$.

If $\text{div}(s)$ is \textit{not} left invariant by the involution
$\sigma$, then consider the meromorphic section of
$\sigma^*\overline{L}$ defined by $\sigma^*s$. (Recall that
$\sigma^*\overline{L}$ and $\sigma^*L$ are identified as real rank two
$C^\infty$ bundles, and the section of $\sigma^*\overline{L}$ defined by
$\sigma^*s$ using this identification is meromorphic.)

Now, fix a holomorphic isomorphism
\begin{equation}\label{sec.-4-def.-alpha}
\alpha\, :\, L\, \longrightarrow\, \sigma^*\overline{L}
\end{equation}
such that the composition
\begin{equation}\label{sec.-4-composition}
L\, \stackrel{\alpha}{\longrightarrow}\,  \sigma^*\overline{L}
\, \stackrel{\sigma^*\overline{\alpha}}{\longrightarrow}\,
\overline{\sigma^*\sigma^*\overline{L}}\,=\, L
\end{equation}
is the identity automorphism of $L$, where $\overline{\alpha}$
is the isomorphism of $\overline{L}$ with
$\overline{\sigma^*\overline{L}}$ induced by $\alpha$. Note that
such an isomorphism exists. Indeed, if
\begin{equation*}
\alpha' \, :\, L\, \longrightarrow\, \sigma^*\overline{L}
\end{equation*}
is any isomorphism, then the automorphism
$\sigma^*\overline{\alpha'}\circ \alpha'$ of $L$ (defined as in
\eqref{sec.-4-composition}) is the multiplication by a nonzero scalar $c\in
{\mathbb C}$. Take any $c_0\in {\mathbb C}$ such that $c^2_0 = c$. Now
the isomorphism $\alpha \,=\, \alpha'/c_0$ satisfies the condition that
the composition in \eqref{sec.-4-composition} is the identity
automorphism of $L$.

Let $s'$ be the meromorphic section of $L$ defined by the above section
$\sigma^*s$ using this isomorphism. Consider the meromorphic section
$s'+s$ of $L$. Since $\text{div}(s)$ is not left invariant by $\sigma$,
this meromorphic section $s'+s$ is not identically zero. The divisor
$\text{div}(s+s')$ is clearly left invariant by the involution $\sigma$.
Hence $L\in J_0(X)$ is in the image of $\rho$. This completes the proof
of the theorem.\hfill \ab
\end{proof}

\section{Nonorientable line bundle}
\setcounter{theo}{0}
In this section we will define a line bundle on $Y$ intrinsically
without using $X$.

Let $\{U_i\}_{i\in I}$ be a covering of $Y$ by open subsets and
for each $U_i$,
\begin{equation*}
\phi_i\, :\, U_i \,\longrightarrow\, {\mathbb R}^2,
\end{equation*}
a $C^\infty$ coordinate chart. Consider the trivial (real) line
bundle $U_i\times {\mathbb R}$ on each $U_i$. Using
\begin{equation*}
\frac{\det \md (\phi_j\circ \phi^{-1}_i)}{\vert\!\det\md (\phi_j\circ
\phi^{-1}_i)\vert}\,\in\, \pm 1\, \subset\, \text{Aut}({\mathbb R})
\end{equation*}
as the transition function over $U_i\cap U_j$ for the pair $(i\, ,j)$,
we get a real line bundle over $Y$. This line bundle will be denoted by
$\xi$. Since the transition functions are $\pm 1$, the line bundle
$\xi^{\otimes 2}$ has a natural isomorphism with the trivial line bundle
$Y\times {\mathbb R}$. Let
\begin{equation}\label{iso.-triv.}
\lambda\, :\, \xi^{\otimes 2}\, \longrightarrow\, Y\times
{\mathbb R}
\end{equation}
be the isomorphism.

We will give a construction of the line bundle $\xi$ without using
coordinate charts. Consider the complement $\bigwedge^2 TY
\setminus\{0_{Y}\}$ of the zero section of the real line bundle $\bigwedge^2
TY$, where $TY$ is the real tangent bundle of $Y$. The multiplicative
group
\begin{equation*}
{\mathbb R}^+\, :=\, \{c\in {\mathbb R}\,\vert\, c >0\}
\end{equation*}
acts on $\bigwedge^2 TY \setminus\{0\}$. The action of any $c\in
{\mathbb R}^*$ sends any $v\in \bigwedge^2 TY \setminus\{0\}$ to $cv$.
Also, the multiplicative group $\pm 1$ acts on $\bigwedge^2 TY
\setminus\{0\}$ by sending any $v$ to $\pm v$. Since these two actions
commute, we have an action of the multiplicative group $\pm 1$ on
\begin{equation*}
Z \, :=\, \frac{\bigwedge^2 TY \setminus\{0_{Y}\}}{{\mathbb R}^+}\, .
\end{equation*}
Now, we have
\begin{equation*}
\xi\, =\, \frac{Z\times {\mathbb R}}{\pm 1},
\end{equation*}
where $\pm 1$ acts diagonally and it acts on ${\mathbb R}$ as
multiplication by $\pm 1$.

We will show that the Klein surface (nonorientable complex) structure on
$Y$ gives an isomorphism of $TY$ with $TY\otimes\xi$, where $TY$ as
before is the (real) tangent bundle of $Y$. To construct the
isomorphism, take a compatible coordinate chart
\begin{equation*}
\phi_i\, :\, U_i \,\longrightarrow\, {\mathbb C}
\end{equation*}
compatible with the nonorientable complex structure. The orientation of
the complex line $\mathbb C$ induces an orientation of $U_i$ using
$\phi_i$. This gives a trivialization of $\xi$ over $U_i$ (this induced
trivialization is also clear from the first construction of $\xi$).
Using $\phi_i$ we have a complex structure on $U_i$ obtained from the
complex structure of $\mathbb C$. Let
\begin{equation*}
\gamma_i\, :\, TU_i \, \longrightarrow\, TU_i\otimes\xi\vert_{U_i}
\end{equation*}
be the isomorphism defined by the almost complex structure of $U_i$ and
the trivialization of $\xi\vert_{U_i}$. If $\phi_j$ is another
compatible coordinate chart then the function $\phi_i\circ\phi^{-1}_j$
is either holomorphic or anti-holomorphic. This immediately implies
that the isomorphism
\begin{equation*}
\gamma_j\, :\, TU_j \, \longrightarrow\, TU_j\otimes\xi\vert_{U_j}
\end{equation*}
(obtained by repeating the construction of $\gamma_i$ for the new
compatible coordinate chart) coincides with $\gamma_i$ over $U_i\cap
U_j$. Consequently, the locally defined isomorphisms $\{\gamma_i\}$
patch together compatibly to give a global isomorphism
\begin{equation}\label{sec.-5-def.-gamma}
\gamma\, :\, TY \, \longrightarrow\, TY\otimes\xi
\end{equation}
over $Y$.

A \textit{nonorientable complex line bundle} over $Y$ is a $C^\infty$
real vector bundle of rank two over $Y$ together with a $C^\infty$
isomorphism of vector bundles
\begin{equation}\label{5.-def.-tau}
\tau \, :\, E\, \longrightarrow\, E\otimes\xi
\end{equation}
satisfying the condition that the composition
\begin{equation}\label{sec.-5-comp.--1}
E\, \stackrel{\tau}{\longrightarrow}\, E\otimes\xi
\, \stackrel{\tau\otimes\text{Id}_\xi}{\longrightarrow}\,
(E\otimes\xi)\otimes\xi\, =\, E\otimes \xi^{\otimes 2}
\, \stackrel{\text{Id}_E\otimes\lambda}{\longrightarrow}\, E
\end{equation}
coincides with the automorphism of $E$ defined by multiplication with
$-1$, where $\lambda$ is defined in \eqref{iso.-triv.}.

Therefore, if for a point $y\in Y$ we fix $w\,\in\, \xi_y$ with $\lambda
(w\otimes w) =1$, then the automorphism of the fiber $E_y$ defined by
\begin{equation*}
v\, \longmapsto\, \langle\tau(v)\, , w^*\rangle
\end{equation*}
is an almost complex structure on $E_y$, where $\langle -\, , -\rangle$
denotes the contraction of $\xi_y$ with its dual line $\xi^*_y$ and $w^*\in
\xi^*_y$ is the dual element of $w$, that is, $\langle w\, , w^*\rangle
=1$.

Let $(E,\tau)$ be a nonorientable complex line bundle over $Y$ as
above. It is easy to see that the $C^\infty$ vector bundle $E$ is
\textit{not} orientable. Indeed, the two orientations on a
two-dimensional real vector space $V$ defined by $J$ and $-J$, where $J$ is
an almost complex structure on $V$, are opposite to each other. To
explain this, note that an orientation of the tangent space $T_yY$, where
$y\in Y$, induces an orientation of the fiber $E_y$ and conversely.
Indeed, giving an orientation of $T_yY$ is equivalent to giving a vector
in $w\,\in\, \xi_y$ with $\lambda (w\otimes w) =1$. As it was shown
above, such an element $w$ gives an almost complex structure on $E_y$.
Hence $E_y$ gets an orientation. Conversely, if we have an orientation
of the fiber $E_y$, then choose the element $w\,\in\, \xi_y$, with
$\lambda (w\otimes w) =1$, that induces this orientation using $\tau$.
Now, $w$ gives an orientation of $T_yY$. Therefore, giving an
orientation of $E_y$ is equivalent to giving an orientation of $T_yY$.
Since the tangent bundle $TY$ is not orientable, we conclude that the
vector bundle $E$ is not orientable.

The total space of the vector bundle $E$ will also be denoted by $E$.
Let
\begin{equation*}
f\, :\, E\, \longrightarrow\, Y
\end{equation*}
be the natural projection. Note that the relative tangent bundle for $f$
(that is, the kernel of the differential $\md f$) is identified with
$f^*E$. So we have the following exact sequence of vector bundle
\begin{equation}\label{5-rel-tan.}
0\, \longrightarrow\, f^*E \, \longrightarrow\, TE \, \longrightarrow
\, f^*TY \, \longrightarrow\, 0
\end{equation}
over the manifold $E$.

The line bundle $f^*\xi$ will be denoted by $\hat{\xi}$. Let
\begin{equation*}
J\, :\, TE\, \longrightarrow\,TE\otimes \hat{\xi}
\end{equation*}
be an isomorphism such that the composition
\begin{equation*}
TE\, \stackrel{J}{\longrightarrow}\, TE\otimes\hat{\xi}
\, \stackrel{J\otimes\text{Id}_{\hat{\xi}}}{\longrightarrow}\,
(TE\otimes\hat{\xi})\otimes\hat{\xi} \, =\,
TE\otimes \hat{\xi}^{\otimes 2}
\, \stackrel{\text{Id}_E\otimes f^*\lambda}{\longrightarrow}\, TE
\end{equation*}
coincides with the automorphism of $E$ defined by multiplication with
$-1$. Assume that the isomorphism $J$ satisfies the following further
conditions:\vspace{.3pc}
\begin{enumerate}
\leftskip .5pc
\item The subbundle $f^*E$ in \eqref{5-rel-tan.} is preserved by $J$
and $J\vert_{f^*E}$ coincides with the isomorphism $f^*\tau$, where
$\tau$ is defined in \eqref{5.-def.-tau}.
\item The action of $J$ on
the quotient $f^*TY$ in \eqref{5-rel-tan.} coincides with the
isomorphism $f^*\gamma$, where $\gamma$ is constructed in
\eqref{sec.-5-def.-gamma}.\vspace{-.4pc}
\end{enumerate}

A holomorphic structure on the nonorientable complex line bundle $E$ is
an isomorphism $J$ as above satisfying the following conditions (apart
from the above conditions) described below.

If we take a coordinate chart $(U,\phi)$ on $Y$ compatible with the
nonorientable Riemann surface structure, then as we saw before, the
restriction $\xi\vert_U$ gets a trivialization. This in turn gives a
trivialization of $\hat{\xi}$ over $f^{-1}(U)$. Using this
trivialization of $\hat{\xi}\vert_{f^{-1}(U)}$, the isomorphism
$J\vert_{f^{-1}(U)}$ becomes an automorphism $J_\phi$ of
$(TE)_{f^{-1}(U)}$ with the property that $J_\phi\circ J_\phi$ coincides
with the automorphism of $(TE)_{f^{-1}(U)}$ given by multiplication with
$-1$. In other words, $J_\phi$ is an almost complex structure on
$f^{-1}(U)$.

A \textit{holomorphic structure} on the nonorientable complex line
bundle $E$ is an isomorphism $J$ satisfying the following two conditions
(apart from the earlier conditions):
\begin{enumerate}
\leftskip .5pc
\item The almost complex structure $J_\phi$ on $f^{-1}(U)$ is
integrable for every compatible coordinate chart.

\item There is a homomorphic isomorphism
\begin{equation*}
f_\phi \, :\, f^{-1}(U) \,\longrightarrow\,
\phi (U)\times{\mathbb C}\, \subset\, {\mathbb C}\times
{\mathbb C}
\end{equation*}
that fits in a commutative diagram
\begin{equation*}
\begin{matrix}
f^{-1}(U) & \stackrel{f_\phi}{\longrightarrow} &
\phi (U)\times{\mathbb C}\\
\Big\downarrow\vcenter{\rlap{$f$}} &&\Big\downarrow\\
U & \stackrel{\phi}{\longrightarrow} & \phi (U)
\end{matrix}
\end{equation*}
(the right vertical arrow is the projection to the first coordinate),
and the restriction of $f_{\phi}$ to any fiber of $f$ is a complex
linear isomorphism with $\mathbb C$.
\end{enumerate}

A \textit{holomorphic line bundle} over $Y$ is defined to be a complex
line bundle equipped with a holomorphic structure.

As in \S\ref{sect.-1}, let $\pi\, :\, X\, \longrightarrow\, Y$ be the
double cover of the nonorientable Riemann surface $Y$ given by local
orientations. As before, let $\sigma$ denote the anti-holomorphic
involution of the Riemann surface $X$.

\begin{theor}[\!]\label{theorem-5.1}
The space of all holomorphic line bundles over $Y$ are in
bijective correspondence with the holomorphic
line bundles $L$ over $X$ with
the property that $\sigma^*\overline{L}$ is holomorphically
isomorphic to $L$.
\end{theor}

\begin{proof}
Let $L$ be a holomorphic line bundle over $X$ such that
$\sigma^*\overline{L}$ is holomorphically isomorphic to $L$. Fix an
isomorphism
\begin{equation*}
\alpha\, :\, L\, \longrightarrow\, \sigma^*\overline{L}
\end{equation*}
as in \eqref{sec.-4-def.-alpha} such that the composition in
\eqref{sec.-4-composition} is the identity automorphism of $L$.

Since the underlying $C^\infty$ line bundle for $\overline{L}$ is
identified with that of $L$, the isomorphism $\alpha$ gives a $C^\infty$
isomorphism of $L$ with $\sigma^*L$ whose composition with itself is the
identity automorphism of $L$. In other words, $\alpha$ is a $C^\infty$
lift to $L$ of the involution $\sigma$ of $X$. Therefore, the quotient
$L/\alpha$ is a real vector bundle of rank two over $X/\sigma \, =\, Y$.
This real vector bundle of rank two over $Y$ will be denoted by $E$.

To construct a complex structure on $E$, first note that the (real) line
bundle $\pi^*\xi$ over $X$ is canonically trivialized, i.e., there is a
natural isomorphism of $\pi^*\xi$ with the trivial line bundle $X\times
{\mathbb R}$ over $X$. Indeed, this follows immediately from the
definitions of $X$ and $\xi$. The complex structure on the fibers on $L$
give an isomorphism
\begin{equation*}
L\, \longrightarrow\, L
\end{equation*}
defined by multiplication by $\sqrt{-1}$. Consider the composition
\begin{equation*}
L\, \longrightarrow\, L \, \longrightarrow\, L\otimes_{\mathbb R}
(X\times {\mathbb R}) \, \longrightarrow\,L\otimes_{\mathbb R}
\pi^*\xi
\end{equation*}
which we denote by $J_0$. Since $\pi^*\xi$ is the pullback of a
line bundle over $Y$, there is a natural lift of the involution $\sigma$
to $\pi^*\xi$. On the other hand, $\alpha$ is a $C^\infty$ lift of the
involution $\sigma$ to $L$. Therefore, we have a lift of the involution
$\sigma$ to $L\otimes_{\mathbb R} \pi^*\xi$. It is straight-forward to
check that the isomorphism $J_0$ defined above commutes with the lifts
of the involution $\sigma$ to $L$ and $L\otimes_{\mathbb R}\pi^*\xi$.
This immediately implies that the isomorphism $J_0$ descends to an
isomorphism of $E$ with $E\otimes_{\mathbb R}\xi$ over $Y$. This
isomorphism of $E$ with $E\otimes_{\mathbb R}\xi$, which we denote
by $J$, clearly satisfies the condition that the composition in
\eqref{sec.-5-comp.--1} is multiplication by $-1$. Therefore, $(E,J)$
is a nonorientable complex line bundle.

It is easy to see that $J$ defines a holomorphic structure on $E$.
Indeed, this is an immediate consequence of the fact that the almost
complex structure on the total space of $L$ is integrable.

For the converse direction, take a holomorphic line bundle $(E,J)$
over $Y$. Consider the (real) rank two $C^\infty$ vector bundle $\pi^*E$
over $X$. Since $\pi^*\xi$ is identified with the trivial line bundle,
the complex structure $\tau$ on $E$ (defined in \eqref{5.-def.-tau})
gives a complex structure on $\pi^*E$. For the same reason, $J$ defines
an integrable complex structure on the total space of $\pi^*E$. Using
the conditions on $J$ the vector bundle $\pi^*E$ gets the structure of a
holomorphic line bundle over $X$.

Since $\pi^*E$ is the pullback of a vector bundle over $Y$, the
involution $\sigma$ of $X$ has a natural $C^\infty$ lift to $\pi^*E$.
The isomorphism of $\pi^*E$ with $\sigma^*\pi^*E$ defined by this lift
gives a holomorphic isomorphism of the holomorphic line bundle $\pi^*E$
with $\sigma^*\overline{\pi^*E}$. This completes the proof of the
theorem.
\end{proof}

\end{document}